%% file: mm_dist_lapl.tex
\documentclass[12pt]{article}
\usepackage{fullpage,graphicx,psfrag,amsmath,amsfonts,verbatim,url}
\usepackage[small,bf]{caption}
\input defs.tex
\usepackage{algorithm,algorithmicx}
\usepackage{xcolor}
\usepackage[noend]{algpseudocode}
\usepackage{float}

\newcounter{algorithmctr}[section]
\renewcommand{\thealgorithmctr}{\thesection.\arabic{algorithmctr}}
\newenvironment{algdesc}%
   {\refstepcounter{algorithmctr}\begin{list}{}{%
       \setlength{\rightmargin}{0\linewidth}%
       \setlength{\leftmargin}{.05\linewidth}}%
       \rmfamily\small
       \item[]{\setlength{\parskip}{0ex}\hrulefill\par%
        \nopagebreak{\bfseries\textsf{Algorithm \thealgorithmctr~}}}}%
   {{\setlength{\parskip}{-1ex}\nopagebreak\par\hrulefill} \end{list}}

\title{Distributed Majorization-Minimization for Laplacian Regularized Problems}
\author{Jonathan Tuck \and David Hallac \and Stephen Boyd}

\begin{document}
\maketitle

\begin{abstract}
We consider the problem of minimizing a block separable convex function 
(possibly nondifferentiable, and including constraints) plus Laplacian 
regularization, a problem that arises in applications including model 
fitting, regularizing stratified models, and multi-period portfolio optimization. 
We develop a distributed majorization-minimization method for this general 
problem, and derive a complete, self-contained, general, and simple proof of 
convergence. Our method is able to scale to very large problems, and we 
illustrate our approach on two applications, demonstrating its scalability and
accuracy.
\end{abstract}

\section{Introduction}\label{s-intro}

Many applications, ranging from multi-period portfolio optimization \cite{BB:17} 
to joint covariance estimation \cite{HPBL:17}, can be modeled as convex optimization
problems with two objective terms, one that is block separable and the other a 
Laplacian regularization term \cite{YGL:16}. The block separable term can be 
nondifferentiable and may include constraints. The Laplacian regularization term is 
quadratic, and penalizes differences between individual variable components. These 
types of problems arise in several domains, including signal processing \cite{PC:17}, 
machine learning \cite{ST:17}, and statistical estimation or data fitting problems with an 
underlying graph prior \cite{AZ:06,MB:11}. As such, there is a need for scalable 
algorithms to efficiently solve these problems.

In this paper we develop a distributed method for minimizing a block-separable 
convex objective with Laplacian regularization.
Our method is iterative; in each iteration a convex problem is solved
for each block, and the variables are then shared with each block's neighbors
in the graph associated with the Laplacian term.
Our method is an instance of a standard and well known general method, 
majorization-minimization (MM) \cite{L:16}, which recovers a wide variety of existing 
methods depending on the choice of majorization \cite{SBP:17}.
In this paper, we derive a diagonal quadratic majorizer of the given Laplacian 
objective term, which has the benefit of separability. This separability
allows for the minimization step in our MM algorithm to be carried 
out in parallel on a block-by-block basis. 
We develop a completely self-contained proof of convergence of 
our method, which relies on no further assumption than the existence of a 
solution.
Finally, we apply our method to two separate applications, multi-period portfolio 
optimization and joint covariance estimation, demonstrating the scalable
performance of our algorithm.

\subsection{Related work}
There has been extensive research on graph Laplacians and Laplacian regularization 
\cite{GR:01,WSQL:07,RHL:13}, and on developing solvers specifically for use in 
optimization over graphs \cite{HWDSBL:15}. In addition, there has been much research 
done on the MM algorithm \cite{AZ:06,RHL:13,L:16,SBP:17}, including
interpreting other well studied algorithms, such as the 
concave-convex procedure and the expectation-maximization algorithm 
\cite{YR:03,WL:11} as special cases of MM.
We are not aware of any previous work that applies the MM 
algorithm to Laplacian regularization.

There has also been much work on the two specific application examples
that we consider.
Multi-period portfolio optimization is studied in, for example,
\cite{AC:00,SB:09,BB:17}, although scalability remains an issue in these
studies.
Our second application example arises in signal processing, specifically
the joint estimation of inverse covariance matrices, which has 
been studied and applied in many different 
contexts, such as cell signaling \cite{FHT:08,DWW:14}, 
statistical learning \cite{BGA:08}, and 
radar signal processing \cite{DBLP:17}. Again, scalability here is either 
not referenced or is still an ongoing issue in these fields.

\subsection{Outline}
In~\S\ref{s-formulation} we set up our notation, and describe the problem of 
Laplacian regularized minimization. 
In~\S\ref{s-majorizers} we show how to construct a diagonal quadratic majorizer of 
a weighted Laplacian quadratic form.
In~\S\ref{s-MM} we describe our distributed MM algorithm, and give a 
complete and self-contained proof of convergence.
Finally, in~\S\ref{s-examples} we present numerical results for 
two applications which demonstrates the effectiveness of our method.

\section{Laplacian regularized minimization}\label{s-formulation}
We consider the problem of minimizing a convex function plus Laplacian 
regularization,
\BEQ \label{e-lrmp}
    \begin{array}{ll}
    \mbox{minimize}   &  F(x) = f(x) + \mathcal{L}(x),
    \end{array}
\EEQ
with variable $x \in \reals^n$.
Here $f: \reals^n \to \reals \cup \{ \infty \}$ is a proper closed 
convex function \cite{R:70,BL:00},
and $\mathcal{L}: \reals^n \rightarrow \reals$ is the Laplacian regularizer
(or Dirichlet energy \cite{E:10}) $\mathcal L(x)=(1/2)x^TLx$,
where $L$ is a weighted Laplacian matrix, \ie, $L=L^T$,
$L_{ij}\leq 0$ for $i \neq j$, and $L\ones =0$, where $\ones$ is the vector
with all entries one \cite{GR:01}.
Associating with $\mathcal L$ the graph with vertices $1,\ldots,n$,
edges indexed by pairs $(i,j)$ with $i<j$ and $L_{ij} \neq 0$, and
(nonnegative) edge weights $w_{ij}=-L_{ij}$,
the Laplacian regularizer can be expressed as
\[
\mathcal L(z) = \sum_{(i,j)\in \mathcal E} w_{ij} (z_i-z_j)^2.
\]

We refer to the problem~(\ref{e-lrmp}) as the \emph{Laplacian regularized 
minimization problem} (LRMP). 
LRMPs are convex optimization problems, which can be solved by a variety of methods,
depending on the specific form of $f$ \cite{BoV:04,NW:06}.
We will let $F^\star$ denote the optimal value of the LRMP.
Convex constraints can be incorporated into LRMP, by defining 
$f$ to take value $+\infty$ when the constraints are violated.
Note in particular that we specifically \emph{do not} assume that the function $f$
is finite, or differentiable (let alone with Lipschitz gradient), or even
that its domain has affine dimension $n$.
In this paper we will make only one additional analytical
assumption about the LMRP~(\ref{e-lrmp}): its sublevel sets 
are bounded. This assumption implies that the LRMP is solvable,
\ie, there exists at least one optimal point $x^\star$, and therefore that its
optimal value $F^\star$ is finite.

A point $x$ is optimal for the LRMP~(\ref{e-lrmp})
if and only if there exists $g \in \reals^n$
such that \cite{R:70,BL:00}
\BEQ \label{e-optimality_condition}
g \in \partial f(x), \quad g + \nabla \mathcal{L}(x) = g + Lx = 0,
\EEQ
where $\partial f(x)$ is the subdifferential of $f$ at $x$ \cite{R:70,C:90}.
For $g \in \partial f(x)$, we refer to 
\[
r=g+Lx
\]
as the \emph{optimality residual} for the LRMP~(\ref{e-lrmp}).
Our goal is to compute an $x$ (and $g\in \partial f(x)$) for which the residual
$r$ is small.

We are interested in the case where $f$ is block separable.
We partition the variable $x$ as $x = (x_1, \ldots, x_p)$, with $x_i \in \reals^{n_i}$,
$n_1+ \cdots + n_p = n$, and assume $f$
has the form 
\[
f(x) = \sum_{i = 1}^p f_i(x_i),
\]
where $f_i: \reals^n \to \reals \cup \{\infty \}$ are closed
convex proper functions.

The main contribution of this paper is a scalable and distributed method 
for solving LRMP in which each of the functions $f_i$ is handled separately.
More specifically, we will see that each iteration of our algorithm
requires the evaluation of a diagonally scaled proximal operator \cite{PB:14}
associated with each block function $f_i$, which can be done in parallel.

\section{Diagonal quadratic majorization of the Laplacian}\label{s-majorizers}
Recall that a function $\hat{\mathcal{L}}: \reals^n \times \reals^n \rightarrow \reals$
is a majorizer of $\mathcal{L}: \reals^n \rightarrow \reals$ 
if for all $x$ and $z$, 
$\hat{\mathcal{L}}(z; z) = \mathcal{L}(z)$, and
$\hat{\mathcal{L}}(x; z) \geq \mathcal{L}(x)$
\cite{L:16,SBP:17}.
In other words, the difference $\hat {\mathcal L}(x,z) - \mathcal L(z)$ is nonnegative,
and zero when $x=z$.

We now show how to construct a quadratic majorizer of the Laplacian regularizer
$\mathcal L$.  This construction is known \cite{SBP:17}, but we give the 
proof for completeness. Suppose $\hat L=\hat L^T$ satisfies 
$\hat{L} \succeq L$, \ie, $\hat L-L$ is positive semidefinite.
The function
\BEQ \label{e-majorizer}
\hat{\mathcal{L}}(x;z) = (1/2)z^T L z + z^T L (x-z) + (1/2)(x-z)^T \hat{L} (x-z),
\EEQ
which is quadratic in $x$, is a majorizer of $\mathcal L$.
To see this, we note that
\[
\begin{array}{ll}
\hat{\mathcal{L}}(x;z) - \mathcal{L}(x) & = (1/2)z^T L z + z^T L (x-z) 
+ (1/2)(x-z)^T \hat{L} (x-z) - (1/2)x^T L x\\
& = (1/2)(x-z)^T(\hat{L} - L)(x-z),
\end{array}
\]
which is always nonnegative, and zero when $x=z$.

In fact, every quadratic majorizer of $\mathcal L$
arises from this construction, for some $\hat L \succeq L$.  To see this we 
note that the difference
$\hat {\mathcal L}(x;z)-\mathcal L(x)$ is a quadratic function of $x$ 
that is nonnegative and zero when $x=z$, so it must have the form
$(1/2)(x-z)^T P(x-z)$ for some $P=P^T\succeq 0$.  
It follows that $\hat {\mathcal L}$ has the form~(\ref{e-majorizer}), with
$\hat L=P+L \succeq L$.

We now give a simple scheme for choosing $\hat L$ in the diagonal 
quadratic majorizer. Suppose $\hat{L}$ is diagonal,
\[
\hat{L} = \diag (\alpha)= \diag(\alpha_1, \ldots, \alpha_n),
\]
where $\alpha \in \reals^n$. 
A simple sufficient condition for
$\hat{L} \succeq L$ is 
$\alpha_i \geq 2L_{ii}$, $i = 1, \ldots, n$.
This follows from standard results for Laplacians \cite{B:98}, 
but it is simple to show directly.
We note that for any $z \in \reals^n$,  we have
\BEAS
z^T (\hat{L} - L) z &=& \sum_{i = 1}^n (\alpha_i - L_{ii}) z_i^2 
+ \sum_{i = 1}^n \sum_{j \neq i}(- L_{ij}) z_i z_j\\
& \geq & \sum_{i = 1}^n L_{ii} z_i^2 + 
\sum_{i = 1}^n \sum_{j \neq i}L_{ij} |z_i| |z_j|\\
&=& |z|^T L |z| \geq 0,
\EEAS
where the absolute value is elementwise.
On the second line we use the inequalities $\alpha_i-L_{ii}\geq L_{ii}$
and for $j \neq i$, 
$-L_{ij}z_iz_j \geq L_{ij}|z_i||z_j|$, which follows since 
$L_{ij} \leq 0$.

In our algorithm described below,
we will require that $\hat L \succ L$, \ie, $\hat L-L$ is positive
definite.  
This can be accomplished by choosing 
\BEQ\label{e-alpha-cond-diags}
\alpha_i > 2L_{ii}, \quad i = 1, \ldots, n. 
\EEQ

There are many other methods for selecting $\alpha$, some of which have additional
properties.
For example, we can choose
$\alpha = 2\lambda_{\max}(L) \ones$,
where $\lambda_{\max}(L)$ denotes the maximum eigenvalue of $L$.  With this choice
we have $\hat L = 2\lambda_{\max}(L) I$.   This diagonal majorization has all
diagonal entries equal, \ie, it is a multiple of the identity.

Another choice (that we will encounter later in \S\ref{ss-lap_reg_est}) takes $\hat L$ to be a 
block diagonal matrix, conformal with the partition of $x$, with each block component
a (possibly different) multiple of the identity, 
\BEQ \label{e-block-scaled-identity-Lhat}
\hat L = \diag( \alpha_1 I_{n_1}, \ldots, \alpha_p I_{n_p}),
\EEQ
where we can take 
\[
\alpha_i > \max_{j \in N_i}  2L_{jj}, \quad i=1, \ldots, p,
\]
where $N_i$ is the index range for block $i$.

\section{Distributed majorization-minimization algorithm}\label{s-MM}

The majorization-minimization (MM) algorithm is an iterative algorithm 
that at each step minimizes a
majorizer of the original function at the current iterate \cite{SBP:17}. 
Since $\hat {\mathcal L}$, as constructed in \S\ref{s-majorizers},
using~(\ref{e-alpha-cond-diags}),
majorizes $\mathcal L$, it follows that
$\hat F = f + \hat {\mathcal L}$ majorizes $F = f+\mathcal L$.
The MM algorithm for minimizing $F$ is then
\BEQ \label{e-update-step}
x^{k+1} = \argmin_{x} \left( f(x) + \hat{\mathcal{L}}(x;x^{k}) \right),
\EEQ 
where the superscripts $k$ and $k+1$ denote the iteration counter.
Note that since $\hat L$ is positive definite, $\hat {\mathcal L}$ is strictly
convex in $x$, so the argmin is unique.

\paragraph{Stopping criterion.}
The optimality condition for the update~(\ref{e-update-step}) is
the existence of $g^{k+1}\in \reals^n$ with
\BEQ \label{e-update-opt}
g^{k+1}\in \partial f(x^{k+1}), \quad g^{k+1}+ \nabla \hat {\mathcal L}
(x^{k+1};x^k )=0.
\EEQ
From $\hat {\mathcal L}(x;z)-\mathcal L(x) = (1/2)(x-z)^T (\hat L-L) (x-z)$, 
we have 
\[
\nabla \hat {\mathcal L}(x^{k+1};x^k) - 
\nabla {\mathcal L}(x^{k+1})  = (\hat L-L) (x^{k+1}-x^k).
\]
Substituting this into~(\ref{e-update-opt}) we get
\BEQ\label{e-lhat-l}
g^{k+1}+ \nabla {\mathcal L} (x^{k+1}) = (\hat L - L)(x^k-x^{k+1}).
\EEQ
Thus we see that 
\[
r^{k+1} = (\hat L-L)(x^k-x^{k+1})
\]
is the optimality residual
for $x^{k+1}$, \ie, the right-hand side of~(\ref{e-optimality_condition}).
We will use $\|r^{k+1}\|_2\leq \epsilon$, where $\epsilon>0$
is a tolerance, as our stopping criterion.
This guarantees that on exit, $x^{k+1}$ satisfies the optimality 
condition~(\ref{e-optimality_condition}) within $\epsilon$.

\paragraph{Absolute and relative tolerance.}
When the algorithm is used to solve problems in which $x^\star$ or $L$
vary widely in size,
the tolerance $\epsilon$ is typically chosen as a combination of an absolute
error $\epsilon_{\mathrm{abs}}$ and a relative error $\epsilon_{\mathrm{rel}}$, 
for example,
\[
\epsilon = \epsilon_{\mathrm{abs}} + 
\epsilon_{\mathrm{rel}} (\|{\hat{L} - L}\|_F + \|x\|_2),
\]
where $\|\cdot\|_F$ denotes the Frobenius norm.

\paragraph{Distributed implementation.}
The update~(\ref{e-update-step}) can be broken down into two steps.
The first step requires multiplying by $L$, and in the other step,
we carry out $p$ independent minimizations in parallel.
We partition the Laplacian matrix $L$ 
into blocks $L_{ij}, i,j = 1, \ldots, p$, conformal with the partition 
$x=(x_1, \ldots, x_p)$.
(In a few places above, we used $L_{ij}$ to denote the $i,j$ entry of
$L$, whereas here we use it to denote the $i,j$ submatrix.  This 
slight abuse of notation should not cause any confusion since
the index ranges, and the dimensions, make it clear whether the entry,
or submatrix, is meant.)
We then observe that our majorizer~(\ref{e-majorizer}) has the form
\[
\hat{\mathcal{L}}(x; z) = 
\sum_{i = 1}^p \hat{\mathcal{L}}_i(x_i; z) + c,
\]
where $c$ does not depend on $x$, and
\[ 
\hat{\mathcal{L}}_i(x_i; z) = (1/2)(x_i - z_i)^T \hat{L}_{ii}(x_i - z_i) 
+ h_i ^T x_i,
\]
where $z_i$ refers to the $i$th subvector of $z$, and 
$h_i$ is the $i$th subvector $L z$,
\[
h_i = L_{ii} z_i + \sum_{j\neq i} L_{ij}z_j.
\]
It follows that 
\[
\hat F(x;x^k) = \sum_{i=1}^p ( f_i(x_i) + \hat {\mathcal L}(x_i;x_i^k) )+c
\]
is block separable.

\begin{algdesc}
        \label{a-distributed-mm}
        \emph{Distributed majorization-minimization.}
        \begin{tabbing}
                {\bf given} Laplacian matrix $L$, and initial
                 starting point $x^0$ in the feasible set of the problem,\\with $f(x^0)<\infty$.\\

                \emph{Form majorizer matrix.} Form diagonal $\hat{L}$ with
                $\hat{L} \succ L$ (using~(\ref{e-alpha-cond-diags})).\\

                {\bf for} $k=1,2, \ldots$ \\

                \qquad \= 1.\ \emph{Compute linear term.} Compute $h^k = L x^k$ 
and residual $r^k = (\hat L - L) (x^{k-1} - x^k)$.\\
                \qquad \= 2. \emph{Update in parallel.} For $i=1,\ldots,p$, update 
each $x_i$ (in parallel) as\\
                        \> \qquad $x_i^{k+1} = \argmin_{x_i} \left( f_i(x_i) + 
                        (1/2)(x_i - x_i^k)^T \hat{L}_{ii}(x_i - x_i^k) + 
                        (h_i^k)^T x_i \right)$.\\
                \qquad \= 3.\ \emph{Test stopping criterion.} Quit if 
                $k \geq 2$ and  $\|r^k\|_2 \leq \epsilon$.
        \end{tabbing}
\end{algdesc}

Step~1 couples the subvectors $x_i^k$; step~2 (the subproblem updates) is 
carried out in parallel for each $i$. We observe that the updates in step~2 
are (diagonally scaled) proximal operator evaluations, \ie, they involve 
minimizing $f_i$ plus a norm squared term, with diagonal quadratic norm; see, 
\eg, \cite{PB:14}. Our algorithm can thus be considered as a distributed 
proximal-based method. 
We also mention that as the algorithm converges (discussed in detail below), 
$x_i^{k+1}-x_i^k \to 0$, which implies that the quadratic terms 
$(1/2)(x_i - x_i^k)^T \hat{L}_{ii}(x_i - x_i^k)$ and their gradients in the 
update asymptotically vanish; roughly speaking, they `go away' as the algorithm 
converges. We will see below, however, that these quadratic terms are critical 
to convergence of the algorithm.

\paragraph{Warm start.} Our algorithm supports \emph{warm 
starting} by choosing the initial point $x^0$ as an estimate of the solution,
for example, the solution of a closely related problem.
Warm starting can decrease
the number of iterations required to converge \cite{YW:02,WB:10};
we will see an example in \S\ref{s-examples}.

\subsection{Convergence}
There are many general convergence results for MM methods, but all of them require 
varying additional assumptions about the objective function \cite{L:16,SBP:17}.
In this section we give a complete self-contained proof of convergence for
our algorithm, that requires no additional assumptions. 
We will show that $F(x^k)-F^\star \to 0$, as $k \to \infty$, and also that
the stopping criterion eventually holds, \ie,
$(\hat L-L)(x^k-x^{k+1})\to 0$.

We first observe a standard result that holds for all MM methods: The objective 
function is non-increasing. We have 
\[
F(x^{k+1}) \leq \hat{F}(x^{k+1};x^{k}) \leq \hat{F}(x^{k};x^{k}) = F(x^{k}),
\]
where the first inequality holds since $\hat F$ majorizes $F$,
and the second since $x^{k+1}$ minimizes $\hat F(x;x^k)$ over $x$.
It follows that $F(x^k)$ converges, and therefore 
$F(x^{k}) - F(x^{k+1}) \to 0$. It also follows that the iterates $x^k$ are bounded, 
since every iterate satisfies $F(x^k)\leq F(x^0)$, and we assume that the sublevel 
sets of $F$ are bounded.

Since $F$ is convex and $g^{k+1} + 
\nabla \mathcal{L}(x^{k+1}) \in \partial F(x^{k+1})$, we have (from the definition
of subgradient)
\[
F(x^k) \geq F(x^{k+1}) + (g^{k+1} + 
\nabla \mathcal{L}(x^{k+1}))^T (x^k - x^{k+1}) .
\]
Using this and~(\ref{e-lhat-l}), we have 
\[
F(x^k) - F(x^{k+1}) \geq (x^k - x^{k+1})^T (\hat{L} - L) (x^k - x^{k+1}).
\] 
Since $F(x^{k})-F(x^{k+1}) \to 0$ as $k\to \infty$,
and $\hat L-L \succ 0$, we conclude that 
$x^{k+1}-x^k \to 0$ as $k\to \infty$.
This implies that our stopping criterion will eventually hold.

Now we show that $F(x^k)\to F^\star$. Let $x^\star$ be any optimal point. Then,
\BEAS
F^\star = F(x^\star) 
&\geq& F(x^{k+1}) + ( (\hat{L} - L) x^{k} - \hat{L}x^{k+1} + L x^{k+1})^T(x^\star-x^{k+1})\\
&=& F(x^{k+1}) + (x^{k} - x^{k+1})^T (\hat L - L) (x^\star-x^{k+1}).
\EEAS
So we have
\[
F(x^{k+1})-F^\star \leq 
-(x^{k} - x^{k+1})^T (\hat L - L) (x^\star-x^{k+1}).
\]

Since $x^{k}-x^{k+1} \to 0$ as $k \to \infty$, and $x^{k+1}$ is bounded,
the right-hand side converges to zero as $k\to \infty$, and so
we conclude $F(x^{k+1})-F^\star \to 0$ as $k\to \infty$.

\subsection{Variations}
\paragraph{Arbitrary convex quadratic regularization.}
While our interest is in the case when $\mathcal L$ is Laplacian regularization,
the algorithm (and convergence proof) 
work when $\mathcal L$ is any convex quadratic, \ie,
$L \succeq 0$, with the choice
\[
\alpha_i > \sum_{j = 1}^n |L_{ij}|, \quad i=1,\ldots, n,
\]
replacing the condition~(\ref{e-alpha-cond-diags}). In fact, the 
condition~(\ref{e-alpha-cond-diags}) is a special case of this condition, for a 
Laplacian matrix.

\paragraph{Nonconvex $f$.} If the objective function in the LRMP is nonconvex, \ie, 
$f$ is nonconvex, then the method proposed in this paper can be extended as a 
heuristic for solving (\ref{e-lrmp}) for nonconvex $f$. It is emphasized that the 
most the algorithm can guarantee is a local optimum, rather than a 
global optimum \cite{BoV:04}.

\section{Examples}\label{s-examples}
In this section we describe two applications of our distributed method for
solving~LRMP, and report numerical results demonstrating 
its convergence and performance.
We run all numerical examples on a 32-core AMD machine with 64 hyperthreads, 
using the Pathos multiprocessing package to carry out computations in parallel
\cite{pathos_multiprocessing}. Our code is available online at 
\url{https://github.com/cvxgrp/mm_dist_lapl}.

\subsection{Multi-period portfolio optimization} \label{ss-MPO}
We consider the problem of multi-period trading with quadratic transaction costs; 
see \cite{BM:14,BB:17} for more detail.
We are to choose a portfolio of $n$ holdings $x_t\in  \reals^n$, for 
periods $t=1,\ldots, T$. We assume the $n$th holding is a riskless holding
(\ie, cash).
We choose the portfolios by solving the problem
\BEQ \label{e-finance-main-prob}
    \begin{array}{ll}
    \mbox{minimize}   &
    \sum_{t = 1}^{T} \left( f_t(x_t)
    + (1/2)(x_t - x_{t-1})^T D_t (x_t - x_{t-1})\right),
    \end{array}
\EEQ 
where $f_t$ is the convex objective function (and constraints) for the portfolio
in period $t$, and the $D_t$'s are diagonal positive definite matrices. 
The initial portfolio $x_0$ is given and constant; $x_1, \ldots, x_T$ 
are the variables. The quadratic term 
$(1/2)(x_t - x_{t-1})^T D_t (x_t - x_{t-1})$ is the transaction cost, 
\ie, the additional cost of trading to move from the previous portfolio 
$x_{t-1}$ to the current one $x_t$.
We will assume that there is no transaction cost associated with cash,
\ie, $(D_t)_{nn}=0$.

The objective function $f_t$ typically includes negative expected return, 
one or more risk constraints or risk avoidance terms, shorting or borrow costs,
and possibly other terms.  It also can include constraints, such as the normalization
$\ones^Tx_t=1$ (in which case $x_t$ are referred to as the portfolio weights),
limits on the holdings or the leverage of the portfolio, or a specified final 
portfolio; see \cite{BM:14,BB:17} for more detail.

We can express the transaction cost as Laplacian regularization on 
$x=(x_1, \ldots, x_T) \in \reals^{Tn}$, plus a quadratic term involving $x_1$,
\begin{eqnarray*}
\lefteqn{\sum_{t = 1}^{T} (1/2)(x_t-x_{t-1})^T D_t (x_t-x_{t-1})}\\
&=&
(1/2) x^T L x + (1/2)x_1^T D_1 x_1 - (D_1 x_0)^T x_1 
+ (1/2) x_0^T D_1 x_0.
\end{eqnarray*}
(Recall that the initial portfolio $x_0$ is given.)
The Laplacian matrix $L$ has block-tridiagonal form given by
\[ L =
\begin{bmatrix}
    D_2 & -D_2 & 0 & \ldots  & 0 & 0 & 0\\
    -D_2 & D_2+D_3 & -D_3 & \ldots & 0 & 0 & 0\\
    0 & -D_3 & D_3+D_4 & \ldots  & 0 & 0 & 0\\
    \vdots & \vdots & \vdots & \ddots & \vdots & \vdots & \vdots\\
    0 & 0 & 0 & \ldots & D_{T-2}+D_{T-1} & -D_{T-1} & 0\\
    0 & 0 & 0 & \ldots & -D_{T-1} & D_{T-1}+D_T & -D_T \\
    0 & 0 & 0 & \ldots & 0  & -D_T & D_T
\end{bmatrix}.
\]
If we assume that the initial portfolio is cash, \ie,
$x_0$ is zero except in the last component, then two of the three extra
terms, $(D_1x_0)^T x_1$ and $(1/2) x_0^T D_1 x_0$, both vanish.
If we lump the extra terms that depend on $x_1$ into $f_1$, the multi-period 
portfolio optimization problem~(\ref{e-finance-main-prob}) has the LRMP form, 
with $p=T$ and $n_i=n$.  The total number of (scalar) variables is $Tn$.
The graph associated with the Laplacian is the simple chain graph; 
roughly speaking, each portfolio $x_t$ is linked to its predecessor 
$x_{t-1}$ and its successor $x_{t+1}$ by the transaction cost.

We can give a simple interpretation for the subproblem update in our method. 
The quadratic term of the subproblem update (which asymptotically goes away as
we approach convergence) adds diagonal risk; the linear term $h_t$ contributes an 
expected return to each asset. These additional risk and return terms come from 
both the preceding and the successor portfolios; they `encourage' the 
portfolios to move towards each other from one time period to the next, 
so as to reduce transaction cost. Each subproblem update
minimizes negative risk-adjusted return, with the given return vector modified
to encourage less trading.

\subsubsection{Problem instance}
We consider a problem with $n=1000$ assets and $T=30$ periods, so the total
number of (scalar) variables is $30000$. The objective functions $f_t$ include a negative
expected return, a quadratic risk given by a factor (diagonal plus low rank) model
with 50 factors \cite{CR:83,BB:17}, and a linear shorting cost. We additionally impose
the normalization constraint $\ones^Tx_t=1$, so the portfolios $x_t$ represent
weights.
The objective functions $f_t$ have the form
\BEQ
f_t(x) = 
-\mu_t^T x + \gamma x^T \Sigma_t x + s_t^T (x)_{-}, \quad t = 1, \ldots, T-1.
\EEQ
Here, $\gamma > 0$ is the \textit{risk aversion parameter},
$\mu_t$ is the expected return, $\Sigma_t$ is the return covariance,
and $s_t$ is the (positive) shorting cost coefficient vector.
The covariance matrices $\Sigma_t$ are diagonal plus a rank~50 (factor)
term,
with zero entries in the last row and column (which correspond to the 
cash asset).
We choose all these coefficients and the diagonal transaction cost
matrices $D_t$ randomly, but with realistic values.
In our problem instance, we choose all of these parameters independent
of $t$, \ie, constant.

We take $f_T$ to be the indicator function for the constraint
$x = e_n$ (\ie, $f_T(x) = 0$ if $x=e_n$, and $\infty$ otherwise),
and the initial portfolio is all cash, $x_0=e_n$. 
So in our multi-period portfolio optimization problem we are planning a 
sequence of portfolios that begin and end in cash.

We can see the interpretation of the subproblem updates in~\S\ref{ss-MPO} by 
looking at 
the subproblem objective functions. 
Assuming we choose the diagonal elements of $\hat L$ to 
be $3 L_{ii}$, we can rewrite the subproblem objective function (at time periods 
$t = 2, \ldots, T-1$ and iteration $k$) as
\[
\begin{array}{ll}
x_t^T (\gamma \Sigma_t + (3/2)(D_t + D_{t+1})) x_t 
- (\mu_t + D_t(2x_t^k - x_{t-1}^k) + D_{t+1}(2x_t^k - x_{t+1}^k))^T x_t + c,
\end{array}
\]
where $c$ is some constant that does not depend on $x_t$.
We see that a diagonal risk term is added,
and the mean return $\mu_t$ is offset by terms that depend on the past, 
current, and future portfolios $x_{t-1}^k$, $x_t^k$, and $x_{t+1}^k$.

\subsubsection{Numerical results}

We first solve the problem instance using CVXPY \cite{DB:16} and solver OSQP 
\cite{SBGBB:17}, which is single-thread. The solve time for this baseline method 
was 120 minutes. 

\begin{figure}
\begin{center}
\includegraphics[width=0.8\textwidth]{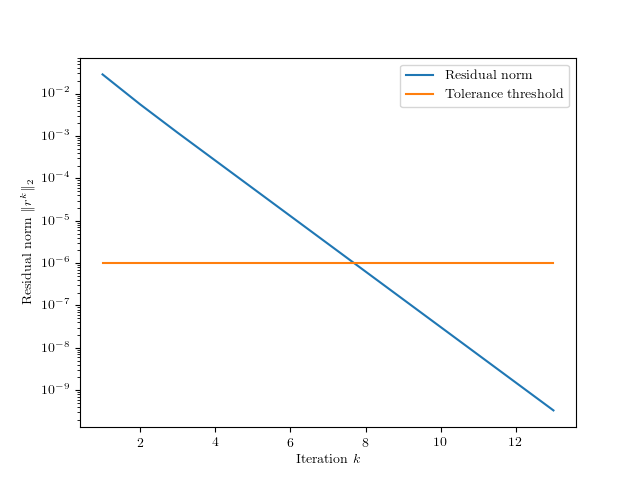}
\end{center}
\caption{Residual norm versus iteration for multi-period portfolio optimization problem.}
\label{f-residual_finance_mm}
\end{figure}

We then solved the problem instance using our method.
We initialized all portfolios as $e_n$, \ie, all cash, and use 
stopping criterion tolerance $\epsilon = 10^{-6}$. 
Our algorithm 
took 8 iterations and 19 seconds to converge, and produced a solution 
that agreed very closely with the CVXPY/OSQP solution. 
Figure~\ref{f-residual_finance_mm} 
shows a plot of the residual norm $\|r^k\|_2$ versus iteration $k$.
Ths plot shows nearly linear convergence, with a reduction in residual
norm by around a factor of $5$ each iteration.

\subsection{Laplacian regularized estimation} \label{ss-lap_reg_est}

We consider estimation of parameters in a statistical model.  
We have a graph, with some data associated with each node; the goal is to fit a model
to the data at each node, with Laplacian regularization used to make neighboring
models similar.

The model parameter at node $i$ is $\theta_i\in \reals^{n_i}$.
The vector of all node parameters is 
$\theta = (\theta_1, \ldots, \theta_p) \in \reals^{n}$, with $n=n_1+\cdots+n_p$.
We choose $\theta$ by minimizing a local loss function and regularizer at 
each node, plus Laplacian regularization:
\[
    \begin{array}{ll}
    \mbox{minimize}   &  \sum_{i = 1}^p f_i(\theta_i) + \mathcal{L}(\theta),
    \end{array}
\]
where $f_i(\theta_i) = \ell_i(\theta_i) + r_i(\theta_i)$, 
where $\ell_i: \reals^{n_i} \to \reals$ is the loss function (for example, the 
negative log-likelihood of $\theta_i$) for the data at node $i$, 
and $r_i: \reals^{n_i} \to \reals$ is a regularizer on the parameter $\theta_i$.
Without the Laplacian term, the problem is separable, and corresponds to fitting 
each parameter separately by minimizing the local loss plus regularizer.
The Laplacian term is an additional regularizer that encourages various
entries in the parameter vectors to be close to each other.

\paragraph{Laplacian regularized covariance estimation.}
We now focus on a more specific case of this general problem, Laplacian regularized 
covariance estimation. At each node, we have some number of samples from 
a zero-mean Gaussian distribution on $\reals^d$,
with covariance matrix $\Sigma_i$, assumed positive definite.
We will estimate the natural parameters (as an exponential family),
the inverse covariance matrices $\theta_i = \Sigma_i^{-1}$.
So here we take the node parameters $\theta_i$ to be symmetric positive definite
$d \times d$ matrices, with $n_i=d(d+1)/2$.
(In the discussion of the general case above, $\theta_i$ is a vector 
in $\reals^{n_i}$; in the rest of this section, $\theta_i$ will denote 
a symmetric $d \times d$ martix.)

The data samples at node $i$ have empirical 
covariance $S_i$ 
(which is not positive definite if there are fewer than $d$ samples).
The negative log-likelihood for node $i$ is (up to a constant
and a positive scale factor)
\[
\ell_i(\theta_i) = \Tr(S_i \theta_i) - \log \det \theta_i.
\] 
We use trace regularization on the parameter,
\[
r_i(\theta_i) = \kappa \Tr(\theta_i),
\]
where $\kappa > 0$ is the local regularization hyperparameter.
We note that we can minimize $f_i(\theta_i) = \ell_i(\theta_i) + 
r_i(\theta_i)$ analytically; the minimizer is 
\[
\theta_i = (S_i+ \kappa I)^{-1}.
\]
(See, \eg, \cite{BGA:08}.)

The Laplacian regularization is used to 
encourage neighboring inverse covariance matrices in the given graph 
to be near each other.
It has the specific form
\[
\mathcal L(\theta_1, \ldots, \theta_p) = \lambda \sum_{(i,j)\in \mathcal E} 
\|\theta_i-\theta_j\|_F^2 = \Tr(\theta^T L \theta),
\]
where the norm is the Frobenius norm, $L$ is the associated weighted Laplacian 
matrix for the graph with vertices $1, \ldots, p$ and edges $\mathcal E$, 
and $\lambda \geq 0$ is a hyperparameter that controls the amount of 
Laplacian regularization. 
When $\lambda =0$, the estimation problem is separable, with analytical
solution
\[
\theta_i=(S_i+\kappa I)^{-1}, \quad i=1, \ldots, p.
\]
For $\lambda \to \infty$,
assuming the graph is connected, 
the estimation problem reduces to finding a single covariance matrix for all
the data, with analytical solution 
\[
\theta_i=(S+p \kappa I)^{-1}, \quad i=1, \ldots, p,
\]
where $S=\sum_{j=1}^p S_j$ is the empirical covariance of all the data together.

We choose the majorizer to be block diagonal with
each block a multiple of the identity, as in~(\ref{e-block-scaled-identity-Lhat}).
The update at each node in our algorithm can be expressed as 
minimizing over $\theta_i$ the function
\[
\Tr((S_i+H_i^k) \theta_i) - \log \det \theta_i
+ \kappa \Tr(\theta_i)
+ (\alpha_i/2) \|\theta_i - \theta_i^k\|_F^2,
\]
where 
\[
H^k= L \theta^k.
\]
This minimization can be carried out analytically. By taking the gradient 
of the subproblem objective function with respect to $\theta_i$ and equating to zero, 
we see that
\[
S_i + H_i^k - \theta_i^{-1} + \kappa I + \alpha_i(\theta_i - \theta_i^k) = 0,
\]
or
\[
\theta_i^{-1} - \alpha_i\theta_i = S_i + H_i^k + \kappa I - \alpha_i\theta_i^k.
\]
This implies that $\theta_i$ and $S_i + H_i^k + \kappa I - \alpha_i\theta_i^k$
share the same eigenvectors \cite{WT:09,DWW:14,HPBL:17}. Let $Q_i \Lambda_i Q_i^T$ 
be the eigenvector decomposion of $S_i + H_i^k + \kappa I - \alpha_i\theta_i^k$.
We find that the eigenvalues of $\theta_i, v_{ij}, j = 1, \ldots, n$, are
\[
v_{ij} = (1/2 \alpha_i)\left(-(\Lambda_{i})_{jj} + 
\sqrt{(\Lambda_{i})_{jj}^2 + 4 \alpha_i}\right).
\]
We have $\theta_i^{k+1} = Q_i V_i Q_i^T$,
where $V_i = \diag(v_{i1}, \ldots, v_{in})$. The computational cost per iteration
is primarily in computing the eigenvector decomposition of 
$S_i + H_i^k + \kappa I - \alpha_i \theta_i^k$, which has order $d^3$.

\subsubsection{Problem instance} \label{ss-cov_prob_instance}
The graph is a $15 \times 15$ grid, with 420 edges, so $p=225$.
The dimension of the data is $d=30$, so each $\theta_i$ is a 
symmetric $30 \times 30$ matrix.
The total number of (scalar) variables in our problem instance is 
$225 \times 30(30+1)/2 = 104625$.

We generate the data for each node as follows. 
First, we choose the four corner 
covariance matrices randomly. The other 221 nodes are given covariance matrices 
using bilinear interpolation from the corner covariance matrices. 
We then generate 
20 independent samples of data from each of the node distributions.
(The samples are in $\reals^{30}$, so the empirical covariance matrices are
singular.)
In our problem instance
we used hyperparameter values $\lambda = .053$ and $\kappa = 0.08$, which were 
chosen to give good estimation performance.

\subsubsection{Numerical results}\label{ss-cov_numerical_results}

The problem instance is too large to reliably solve using CVXPY and the solver 
SCS~\cite{ocpb:16}, which stops after two hours with the status 
message that the computed solution may be inaccurate.

\begin{figure}
\begin{center}
\includegraphics[width=0.8\textwidth]{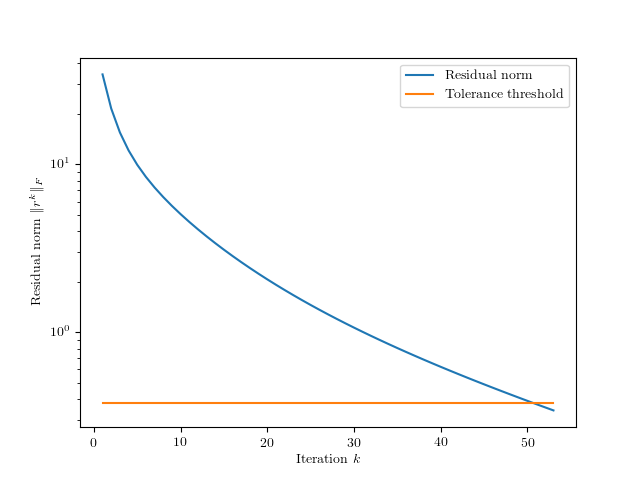}
\end{center}
\caption{Residual norm vs.\ iteration for Laplacian regularized covariance estimation 
problem.}
\label{f-residual_covariance_mm}
\end{figure}

We solved the problem using our distributed method, with absolute tolerance 
$\epsilon_{\mathrm{abs}} = 10^{-5}$ and
relative tolerance $\epsilon_{\mathrm{rel}} = 10^{-3}$.
The method took 54 iterations and 13 seconds to converge. 
Figure~\ref{f-residual_covariance_mm} is a plot of the residual norm $\|r^k\|_F$ 
versus iteration $k$.

\paragraph{Regularization path via warm-start.}
To illustrate the advantage of warm-starting our algorithm,
we compute the entire regularization path, \ie, the solutions of the 
problem for 100 values of $\lambda$, spaced logarithmically between
$10^{-5}$ and $10^4$.

Computing these 100 estimates by running the algorithm for each value of 
$\lambda$ sequentially, without warm-start, requires 26000 total iterations (an 
average of 260 iterations per choice of $\lambda$) and 81 minutes.
Computing these 100 estimates by running the algorithm using warm-start,
starting from $\lambda = 10^{-5}$, requires only 2000 total iterations (an 
average of 20 iterations per choice of $\lambda$) and 7.1 minutes. 
For the specific instance solved above, 
the algorithm converges in
only 2.5 seconds and 10 iterations using warm-start, compared to 
13 seconds and 54 iterations using cold-start.

\begin{figure}
\begin{center}
\includegraphics[width=0.8\textwidth]{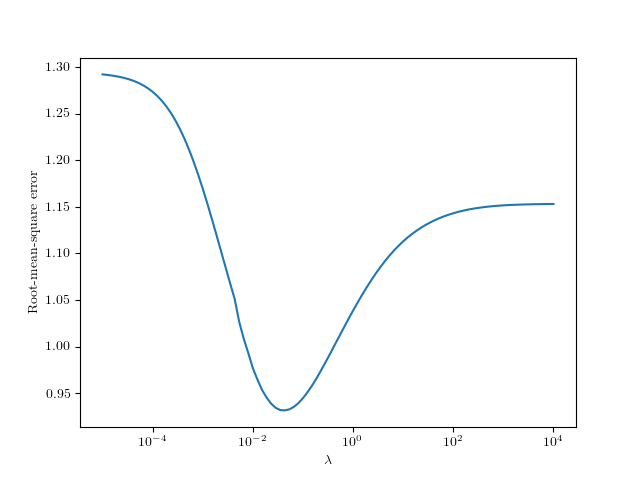}
\end{center}
\caption{Root-mean-square error of the optimal estimates vs. $\lambda$.}
\label{f-covariance_reg_path}
\end{figure}

While the point of this example is the algorithm that computes the estimates, 
we also explore the performance of the method.
For each of the 100 values of $\lambda$ we compute the root-mean-square error 
between our estimate of the inverse covariance and the true inverse covariance,
which we know, since we generated them.
Figure~\ref{f-covariance_reg_path} shows a plot of the root-mean-square error of
our estimate versus the value of $\lambda$.
This plot shows that the method works, \ie, produces better estimates of 
the inverse covariance matrices than handling them separately (small $\lambda$)
or fitting one inverse covariance matrix for all nodes (large $\lambda$).

\section*{Acknowledgements}
The authors would like to thank Peter Stoica for his insights and comments on 
early drafts of this paper. We would also like to thank the Air Force Research
Laboratory, and in particular Muralidhar Rangaswamy, for discussions of
covariance estimation arising in radar signal processing.

\bibliography{template}

\end{document}

%% file: defs.tex
\newcommand{\ones}{\mathbf 1}
\newcommand{\reals}{{\mbox{\bf R}}}

\newcommand{\Tr}{\mathop{\bf Tr}}
\newcommand{\diag}{\mathop{\bf diag}}

\newcommand{\argmin}{\mathop{\rm argmin}}

\newcommand{\eg}{{\it e.g.}}
\newcommand{\ie}{{\it i.e.}}

\newcommand{\BEAS}{\begin{eqnarray*}}
\newcommand{\EEAS}{\end{eqnarray*}}
\newcommand{\BEA}{\begin{eqnarray}}
\newcommand{\EEA}{\end{eqnarray}}
\newcommand{\BEQ}{\begin{equation}}
\newcommand{\EEQ}{\end{equation}}
\newcommand{\BIT}{\begin{itemize}}
\newcommand{\EIT}{\end{itemize}}